\newcommand{\M}{{\mathcal M}} % la varieta' Lorentziana
\newcommand{\R}{I\!\!R}   % i numeri reali
\newcommand{\N}{I\!\!N}   % i numeri naturali
\newcommand{\ddso}{{\frac{\rm d}{{\rm d}s}}\big\vert_{s=0}}
\newcommand{\ind}{{\rm ind}}

\documentclass[oneside,draft,11pt]{amsart}

\usepackage{amsmath}               % AmSLaTeX
\usepackage{amsthm}                % aggiunge nuovi ambienti tipo teorema
\usepackage{times}

\title[The Morse Index Theorem in semi-Riemannian Geometry]%
{A Note on the Morse Index Theorem for Geodesics between
Submanifolds in semi-Riemannian Geometry}
\author[P.\ Piccione]{Paolo Piccione%
\thanks{Partially sponsored by CNPq, Brazil, Processo n.\ 301410/95-0 (RN)}}
\address{Departamento de Matem\'atica, Universidade de S\~ao Paulo, Brazil}
\email{piccione@ime.usp.br}
\urladdr{http://www.ime.usp.br/\~{}piccione}
\author[D.\ Tausk]{Daniel V.\ Tausk}
\address{Departamento de Matem\'atica, Universidade de S\~ao Paulo, Brazil}
\email{tausk@ime.usp.br}

\begin{document}

% Theorems and such

% Section 2

\theoremstyle{plain}\newtheorem{kernlema}{Lemma}[section]
\theoremstyle{plain}\newtheorem{lema}[kernlema]{Lemma}
\theoremstyle{plain}\newtheorem{corJ}[kernlema]{Corollary}
\theoremstyle{definition}\newtheorem{defpartnorm}[kernlema]{Definition}
\theoremstyle{plain}\newtheorem{indexth}[kernlema]{Theorem} 
\theoremstyle{remark}\newtheorem{light}[kernlema]{Remark}
\theoremstyle{plain}\newtheorem{twosubmanifolds}[kernlema]{Theorem}
\theoremstyle{remark}\newtheorem{h1}[kernlema]{Remark}
\theoremstyle{remark}\newtheorem{oldresults}[kernlema]{Remark}
\theoremstyle{remark}\newtheorem{stationary}[kernlema]{Remark}
%%%%%%%%%%%%%%%%%%%%%%%%%%%%%%%%%%%%%%%%%%%%%%%%%%%%%%%%%%%%%%%%%%%%%%%%%%%
%%%%%%%%%%%%%%%%%%%%%%%%%%%%%%%%%%%%%%%%%%%%%%%%%%%%%%%%%%%%%%%%%%%%%%%%%%%

\begin{abstract} 
The computation of the index of the Hessian of the action functional
in semi-Riemann\-ian geometry at geodesics with
two variable endpoints is reduced to the case
of a {\em fixed\/} final endpoint. Using this observation,
we give an elementary  proof of the Morse Index 
Theorem for Riemannian geodesics with two variable endpoints,
in the spirit of the original Morse's proof. This approach
reduces substantially the effort required in the proofs
of the Theorem given in~\cite{A,Bol,Kal}. 
Exactly the same  argument works also in the case of 
timelike geodesics between two submanifolds of a Lorentzian 
manifold. For the extension to the lightlike Lorentzian case, 
just minor changes are required and 
one obtains easily a proof of the focal index theorem
presented in \cite{EK}.
\end{abstract}

\maketitle

%%%%%%%%%%%%%%%%%%%%%%%%%%%%%%%%%%%%%%%%%%%%%%%%%%%%%%%%%%%%%%%%%%%%%%%%%%%
%%%%%%%%%%%%%%%%%%%%%%%%%%%%%%%%%%%%%%%%%%%%%%%%%%%%%%%%%%%%%%%%%%%%%%%%%%%
\begin{section}{Introduction}
\label{sec:intro}
\noindent
A geodesic in a semi-Riemannian manifold $(\M,g)$ is a smooth curve
$\gamma:[a,b]\longmapsto\M$ that is a stationary point for the action functional
$f(z)=\frac12\int_a^b g(\dot z,\dot z)\;{\rm d}t$ defined in the set
of paths $z$ joining two given points of $\M$.
If $(\M,g)$ is Riemannian, i.e., if $g$ is positive definite,
given one such critical point $\gamma$, the celebrated Morse 
Index Theorem relates  some analytical properties of the second variation of 
$f$ at $\gamma$ with some geometrical properties of $\gamma$. 
More precisely, the {\em index\/} of ${\rm Hess}_f$ at $\gamma$,
that gives the number of {\em essentially
different\/} directions in which $\gamma$ can be deformed to obtain a shorter
curve, equals the number of conjugate points along $\gamma$ counted with 
multiplicity, excluding the endpoints $\gamma(a)$ and $\gamma(b)$.

The Index Theorem opened a very active field of research for
both geometers and analysts, and the original result of Morse
was successively extended in several directions.  
Beem and Ehrlich extended the results to the case of timelike
Lorentzian geodesics (see~\cite{BEE}) and to the lightlike
Lorentzian case (\cite{BE, BEE}). The case of a Riemannian
geodesic with endpoints variable in two submanifolds of $\M$ has
been treated by several authors, including
Ambrose, Bolton and Kalish, (see~\cite{A,Bol,Kal}, see also~\cite{Taka}).
Following the approach of Kalish~\cite{Kal},
Ehrlich and Kim have then proven in \cite{EK} the Morse Index Theorem for lightlike
geodesics with endpoints varying on two spacelike submanifolds of 
a Lorentzian manifold. 
The case of spacelike geodesics in semi-Riemannian manifolds was treated
by Helfer in~\cite{Hel}, where an extension of the Index Theorem
was proven in terms of the {\em Maslov\/} index of a curve,
and by the introduction of a notion of {\em signature\/} for conjugate points.
Edwards extended in~\cite{Ed} the Morse Index Theorem 
to the case of formally self-adjoint linear systems of ODE's,
and Smale proved in \cite{Sm} a general version 
of the Index Theorem for strongly elliptic operators on a Riemannian manifold.

The key point in the original Morse's proof of the theorem was
the introduction of a function $i:[a,b]\longmapsto\N$ that gives
the index of the form $I_t$, which is the
Hessian ${\rm Hess}_f$ restricted to the geodesic $\gamma\vert_{[a,t]}$. 
Using a suitable subdivision of the interval $[a,b]$
and some geometrical arguments (see
\cite{M,dC}) Morse proved that $i$ is non decreasing and
left continuous, with discontinuities precisely at the conjugate points,  
and  that the jump of $i$ at each discontinuity point $t_0$ is given by the value
of the multiplicity of the conjugate point $\gamma(t_0)$.

When passing to the case of variable endpoints, i.e., when one admits 
variations with curves having endpoints varying on two fixed submanifolds
$P$ and $Q$ of $\M$, in which case a stationary point of $f$ is a geodesic $\gamma$
that is orthogonal to $P$ and $Q$ at its endpoints,
some obstructions to the use of
the original argument of Morse arise, due mainly to the fact that the 
restricted index form $I_t$ does not detect the influence of the final manifold $Q$. 

Ambrose~\cite{A} gave a proof of the Index Theorem that uses the subdivision argument,
by introducing a family $Q_t$ of {\em localized end-manifolds\/} along $\gamma$,
constructed with the help of the geodesic flow of the normal bundle of $P$ around
$(\gamma(a),\dot\gamma(a))$. This construction leads to technical difficulties (see also
\cite{Taka}),
due to the fact that the submanifold $Q_t$ may lose dimension and differentiability.
The proof of Bolton~\cite{Bol} also uses a subdivision argument, and it
avoids the introduction of the manifolds $Q_t$, but it employs a
restricted index function which is no longer nondecreasing. 

The passage to a restricted index function is avoided in
Kalish's proof of the Index Theorem in the variable endpoints case 
(see~\cite{Kal}).
In this article, it is given an explicit  direct sum decomposition
of the space ${\mathcal H}^{(P,Q)}=B\oplus B^c_+\oplus B^c_-$
of vector fields along $\gamma$ which are everywhere orthogonal to $\gamma$ 
and tangent to $P$ and $Q$ respectively at $\gamma(a)$ and $\gamma(b)$.
The index theorem is deduced with a study of the sign of the index form in
each of the three spaces; the definition of such decomposition is not 
very natural, and the remaining calculations are  rather involved. 

Ehrlich and Kim~\cite{EK} have adapted Kalish's proof to the case of lightlike
Lorentz\-ian geodesics, where a suitable quotient space is used, in analogy
with the null Morse Index Theorem of \cite{BE,BEE}.

The aim of this paper is to show that the proof  of the Morse
Index Theorem for geodesics with two variable endpoints is a simple
adaptation of the classical proof for the fixed endpoints case, 
in the spirit  of the original proof of Morse, which is well understood. 
To this goal, the key observation
is that the case of a geodesic with final point varying on a submanifold
$Q$ can be deduced immediately from the case of a fixed final endpoint
(see Theorem~\ref{thm:twosubmanifolds}) by using a natural splitting
of the space ${\mathcal H}^{(P,Q)}$. 
Moreover, we emphasize that the case of causal (nonspacelike) Lorentzian
geodesics is essentially analogous to the Riemannian case.

We try to keep all the statements and proofs of the paper
at the maximum level of generality; in particular,
we present an approach that unifies
the Riemannian and the causal Lorentzian case, obtaining  
a proof of all the results for Riemannian and causal Lorentzian
geodesics at the same time. In Remark~\ref{thm:oldresults}, 
among other things we  observe that, in the Lorentzian lightlike
case, the use of the quotient bundle employed in \cite{BE, BEE, EK}
is not really essential for the computation of the (non augmented) index, 
which allows to give an easier statement of the focal index theorem. 

It is also important to observe that the result of Theorem~\ref{thm:twosubmanifolds}
applies to a great number of situations in semi-Riemannian geometry
where the Morse Index Theorem may {\em not\/} work, like for instance 
in the case of spacelike geodesics in stationary Lorentzian manifolds 
(see Remark~\ref{thm:stationary}).

\end{section}
%%%%%%%%%%%%%%%%%%%%%%%%%%%%%%%%%%%%%%%%%%%%%%%%%%%%%%%%%%%%%%%%%%%%%%%%%%%
%%%%%%%%%%%%%%%%%%%%%%%%%%%%%%%%%%%%%%%%%%%%%%%%%%%%%%%%%%%%%%%%%%%%%%%%%%%
\begin{section}{The Index Theorem}\label{sec:morse}
Let $(\M,g)$ be a semi-Riemannian manifold, $m={\rm dim}(\M)$,
$P\subset\M$ be a smooth submanifold of $\M$ and $\gamma:[a,b]\longmapsto\M$
be a non constant geodesic in $\M$, with $\gamma(a)\in P$ and $\dot\gamma(a)\in 
T_{\gamma(a)}P^\perp$.  
We will say that $\gamma$ is spacelike, timelike or lightlike according
to $g(\dot\gamma,\dot\gamma)$ positive, negative or zero, respectively;
by {\em causal\/} we will mean either timelike or lightlike.

Let $\nabla$ denote the Levi--Civita connection
of $g$ and let \[R(X,Y)=\nabla_X\nabla_Y-\nabla_Y\nabla_X-\nabla_{[X,Y]}\]
be the curvature tensor of $\nabla$;  moreover, for all $p\in P$ and all 
$n\in T_pP^\perp$, let ${\mathcal S}^P_n$ be the second fundamental 
form of $P$ in the orthogonal direction $n$, which is the following 
symmetric bilinear form on $T_pP$:
\[{\mathcal S}^P_n(v_1,v_2)=g(n,\nabla_{v_1}V_2),\]
where $V_2$ is any extension of $v_2$ to a vector field tangent to $P$.
Observe that we are {\em not\/} in principle making any non degeneracy 
assumption on $P$, but if the metric is non degenerate on $T_pP$ 
then we can also define a linear map
${\mathcal S}^P_n:T_pP\longmapsto T_pP$ such that 
$g({\mathcal S}^P_n(v_1),v_2)={\mathcal S}^P_n(v_1,v_2)$.

Given a (piecewise) smooth vector field $V$ along $\gamma$, we
denote by $V'$ the covariant derivative of $V$ along $\gamma$;
if $V$ is piecewise smooth and $\tau\in[a,b]$, 
the symbols $V'(\tau^-)$ and $V'(\tau^+)$ will mean
respectively the left and right limits of $V'(t)$ as $t\to \tau$.

If $(\M,g)$ is Lorentzian, i.e., if the index of $g$ is $1$,
and $\gamma$ is timelike, we have that
$T_{\gamma(a)}P$ is  spacelike, in the sense that the 
restriction of $g$ to $T_{\gamma(a)}P$
is positive definite. More in general, the restriction of the metric
$g$ to the orthogonal space $\dot\gamma(t)^\perp$ is positive definite
for all $t\in[a,b]$. If $\gamma$ is lightlike, the restriction of 
the metric to the orthogonal space is just positive semi-definite 
(having a one dimensional kernel spanned
by $\dot\gamma(t)$). However, if one assumes that 
$\dot\gamma(a)\not\in T_{\gamma(a)}P$, then again
$T_{\gamma(a)}P$ is spacelike.

Let $\tilde{\mathcal H}^P$ denote the vector space of all
piecewise smooth vector fields $V$ along $\gamma$ such that 
$V(a)\in T_{\gamma(a)}P$ and let ${\mathcal H}^P$
be the subspace of $\tilde{\mathcal H}^P$ consisting of those $V$
such that $g(V,\dot\gamma)\equiv0$ and $V(b)=0$; moreover, let 
$I^P:\tilde{\mathcal H}^P\times\tilde{\mathcal H}^P\longmapsto\R$ 
be the symmetric bilinear form given by:
\begin{equation}\label{eq:bilform}
I^P(V,W)=\int_a^b\Big[g(V',W')+g(R(\dot\gamma,V)\,\dot\gamma,W)\Big]\;{\rm d}t-
{\mathcal S}^P_{\dot\gamma(a)}(V(a),W(a)).
\end{equation}
Observe that if the submanifold $P$ consists of just one point, 
the term involving  its second fundamental
form ${\mathcal S}^P_{\dot\gamma(a)}$ in (\ref{eq:bilform}) 
disappears. In this case we'll write just $I$
instead of $I^P$.

Integration by parts on $g(V',W')$ gives yet another expression for $I^P$:
\begin{equation}\label{eq:bilformint}
\begin{split}
I^P(V,W)=&\int_a^b g(R(\dot\gamma,V)\,\dot\gamma-V'',W) 
\;{\rm d}t+\\&+g(V'(b),W(b))-g(V'(a),W(a))-{\mathcal S}^P_{\dot\gamma(a)}(V(a),W(a))+
\\&+\sum_{i=1}^{N-1}
g(V'(t^-_i)-V'(t^+_i),W(t_i)),
\end{split}
\end{equation}
where $a=t_0<t_1<\ldots<t_N=b$ is a partition of $[a,b]$ such that $V$ is smooth 
in each interval $[t_i,t_{i+1}]$, $i=0,1,\ldots,N-1$.

It is well known that $\gamma$ is a stationary point for the
action functional \[f(z)=\frac12\int_a^bg(\dot z,\dot z)\;{\rm d}t\] defined
in the set $\Omega_{P,\gamma(b)}$ of all piecewise smooth curves 
$z:[a,b]\to\M$ joining
$P$ and $\gamma(b)$. Under the viewpoint of Calculus of Variations and
Global Analysis, the vector space ${\mathcal H}^P$ is a subspace of the
{\em tangent space\/} of $\Omega_{P,\gamma(b)}$ at $\gamma$, and 
$I^P\big\vert_{{\mathcal H}^P}$
is the bilinear form given by the {\em second variation\/} of $f$ at the stationary point
$\gamma$. 
We will be concerned with the {\em index\/} of $I^P$ in ${\mathcal H}^P$, defined as
follows. If ${\mathcal K}$ is a vector subspace of $\tilde{\mathcal H}^P$, then the index
$i(I^P,{\mathcal K})$ of $I^P$ in $\mathcal K$ is the number:
\[\ind(I^P,{\mathcal K})=\sup\{{\rm dim}({\mathcal V}):\ {\mathcal V}\ \text{subspace of}\ 
{\mathcal K}\ \text{with}\ I^P\big\vert_{\mathcal V}<0\},\]
and we set
\begin{equation}\label{eq:MorseInd} 
\ind(I^P)=\ind(I^P,{\mathcal H}^P).
\end{equation}
The number $\ind(I^P)$ will be called the {\em Morse Index\/} of $\gamma$.

A Jacobi field along $\gamma$ is a smooth vector field
$J$ satisfying the linear equation 
$J''-R(\dot\gamma,J)\,\dot\gamma=0$. We say that $J$ is a $P$-Jacobi field if 
it satisfies in addition:
\begin{equation}\label{eq:forgot}
J(a)\in T_{\gamma(a)}P,
\end{equation}
and
\begin{equation}\label{eq:J1}
g(J'(a),w)+{\mathcal S}^P_{\dot\gamma(a)}(J(a),w)=0,\quad 
\text{for all}\ w\in T_{\gamma(a)}P.
\end{equation}
If the metric is non degenerate on $T_{\gamma(a)}P$ we can rewrite 
(\ref{eq:J1}) as \[J'(a)+{\mathcal S}^P_{\dot\gamma(a)}(J(a))\in 
T_{\gamma(a)}P^\perp.\] In this case, a simple counting
argument shows that the dimension of the vector space of $P$-Jacobi 
fields along $\gamma$ is precisely equal to $m$ and that
the dimension of $P$-Jacobi fields satisfying $g(J,\dot\gamma)=0$ 
is equal to $m-1$ (for $P$-Jacobi fields the condition
$g(J,\dot\gamma)=0$ is equivalent to $g(J'(a),\dot\gamma(a))=0$).
Observe that if $P$ is a point, then a $P$-Jacobi field
is simply a Jacobi field $J$ along $\gamma$ such that $J(a)=0$.

Two points $q_0=\gamma(t_0)$ and $q_1=\gamma(t_1)$, $t_0,t_1\in[a,b]$, 
are said to be {\em conjugate\/} along $\gamma$ if there exists a non null
Jacobi field $J$ along $\gamma$ with $J(t_0)=0$ and $J(t_1)=0$.
A point $q_0=\gamma(t_0)$, $t_0\in\,]a,b]$ is said to be a {\em $P$-focal point\/} 
along $\gamma$
if there exists a non null $P$-Jacobi field $J$ along $\gamma$ such that
$J(t_0)=0$; the {\em geometrical multiplicity\/} $\mu^P(t_0)$ of a $P$-focal point
$\gamma(t_0)$ is the dimension of the vector space of all $P$-Jacobi fields
along $\gamma$ that vanish at $t_0$. If $\gamma(t_0)$ is not $P$-focal, we set
$\mu^P(t_0)=0$.

It is well known that,
if $\gamma$ is either Riemannian or causal Lorentzian,
and $\dot\gamma(a)\in T_{\gamma(a)}P^\perp\setminus T_{\gamma(a)}P$
(see Remark~\ref{thm:light}),
then the set of
$P$-focal points along $\gamma$ is discrete,
\footnote{
As proved in~\cite{Hel}, along a spacelike Lorentzian geodesic,
or more in general along a semi-Riemannian geodesic, the conjugate points
may accumulate.}
 hence finite.
Namely, if $J_1,\ldots,J_m$ is a linear basis for the space of $P$-Jacobi fields
along $\gamma$ and
$E_1,\ldots,E_m$ is a parallely transported orthogonal basis along $\gamma$,
then the smooth function $r(t)={\rm det}(g(J_i,E_j))$ has only simple zeroes
on $[a,b]$, i.e., zeroes of finite multiplicity, exactly at those points
$t_0\in[a,b]$ such that $\gamma(t_0)$ is a $P$-focal point along $\gamma$
(see for instance~\cite[Ex.\ 8, p.\ 299]{ON}).
Similarly, for all $q_0=\gamma(t_0)$, the set of points $q_1$ that
are conjugate to $q_0$ along $\gamma$ is finite.  

We are interested in 
the kernel of the restriction of $I^P$ to ${\mathcal H}^P$. To this aim, 
we introduce the spaces  ${\mathcal N}$ and
${\mathcal J}_0$ as follows:
\begin{equation}\label{eq:defN}
\begin{split}
&{\mathcal N}=\Big\{f\dot\gamma:f:[a,b]\longmapsto\R\  \text{piecewise smooth},\ f(a)=f(b)=0\Big\};
\\ \\
&{\mathcal J}_0=\Big\{P\text{-Jacobi fields}\ J\ \text{along}\ \gamma:J(b)=0\Big\}.
\end{split}
\end{equation}
If $\gamma$ is
lightlike we have ${\mathcal N}\subset{\mathcal H}^P$ and in fact ${\mathcal N}$ 
is contained in the kernel of $I^P$ in
${\mathcal H}^P$, as follows directly from (\ref{eq:bilform}). 
We now compute this kernel in the case of Riemannian
or causal Lorentzian  geodesics.

\begin{kernlema}\label{thm:kernlema}
Let $(\M,g)$ be either Riemannian or Lorentzian; in the latter case
assume that $\gamma$ is causal.
The kernel of the restriction of the bilinear form $I^P$ to ${\mathcal H}^P$ 
is equal to ${\mathcal J}_0$ if $(\M,g)$ is Riemannian or if $(\M,g)$ is Lorentzian 
and $\gamma$ is timelike. If $\gamma$ is lightlike and $\dot\gamma(a)\in T_{\gamma(a)}P^\perp
\setminus T_{\gamma(a)}P$, 
this kernel is equal to ${\mathcal J}_0\oplus{\mathcal N}$.
\end{kernlema}

\begin{proof}
Observe that a $P$-Jacobi field which vanishes at some instant on $]a,b]$ is automatically 
orthogonal to $\gamma$, so that we really have ${\mathcal J}_0\subset{\mathcal H}^P$. 
If $V\in{\mathcal H}^P$ is in the kernel
of (the restriction of) $I^P$, it follows from (\ref{eq:bilformint}) 
and usual techniques of calculus of variations that
$V''-R(\dot\gamma,V)\,\dot\gamma$ is parallel to $\dot\gamma$ and that 
$V$ satisfies equation (\ref{eq:J1}). Since
$V''-R(\dot\gamma,V)\,\dot\gamma$ is also orthogonal to $\dot\gamma$, 
it follows that $V$ is a Jacobi field, except for the
case where $\gamma$ is lightlike. In the latter case, 
we get $V''-R(\dot\gamma,V)\,\dot\gamma=\varphi\dot\gamma$
for some function $\varphi$ and therefore $V-f\dot\gamma$ is a 
Jacobi field, where $f$ satisfies $f''=\varphi$ and $f(a)=f(b)=0$.
Observe that ${\mathcal J}_0\cap{\mathcal N}=\{0\}$ 
because $\dot\gamma(a)\not\in T_{\gamma(a)}P$.
\end{proof}

The proof of the Index Theorem for Riemannian or causal Lorentzian
geodesics with initial endpoint varying on a submanifold and fixed 
endpoint is a simple adaptation of the classical Morse proof of the 
Index Theorem in the case of fixed endpoints (see for instance~\cite{dC,M}). 
For the reader's convenience, we outline briefly such adaptation.

We start with the following:
\begin{lema}\label{thm:lema}
Let $J_1,J_2,\ldots,J_n$ be any family of $P$-Jacobi fields (not necessarily
linearly independent) and
$\phi_1,\ldots,\phi_n,\psi_1,\ldots,\psi_n$ be real piecewise smooth functions
on $[a,b]$. Then,
\begin{equation}\label{eq:tesi}
\begin{split}
I^P(\sum_{i=1}^n\phi_i\cdot J_i,\sum_{j=1}^n\psi_j\cdot J_j)=&
\int_a^bg(\sum_{i=1}^n\phi_i'\cdot J_i,\sum_{j=1}^n\psi_j'\cdot J_j)\;{\rm d}t
+\\&+g(\sum_{i=1}^n\phi_i(b)\cdot J_i'(b),\sum_{j=1}^n\psi_j(b)\cdot J_j(b)).
\end{split}
\end{equation}
\end{lema}
\begin{proof}
It is a simple computation that uses the Jacobi equation,
formulas (\ref{eq:J1}), (\ref{eq:bilform}) and the fact that, for $P$-Jacobi fields
$J_i$ and $J_j$, one has $g(J_i',J_j)=g(J_i,J_j')$.
\end{proof}
For Riemannian or causal Lorentzian geodesics, the above Lemma
gives immediately the following Corollary:
\begin{corJ}\label{thm:corJ}
Let $(\M,g)$ be either Riemannian or Lorentzian; in the latter case assume that
$\gamma$ is causal and that $\dot\gamma(a)\in T_{\gamma(a)}P^\perp\setminus
T_{\gamma(a)}P$. Suppose there are no
$P$-focal points along
$\gamma$. Let
$V,J\in\tilde{\mathcal H}^P$ be vector fields orthogonal to
$\gamma$, with $J$ a $P$-Jacobi field and
such that $V(b)=J(b)$. Then $I^P(V,V) \ge I^P(J,J)$. In the Riemannian and timelike Lorentzian case equality holds if and only if $V=J$, and in the lightlike Lorentzian case it holds if and only if $V-J\in{\mathcal N}$.
\end{corJ}
\begin{proof}

Set $k={\rm dim}(P)$. For $i=1,\ldots,k$, choose Jacobi fields $J_i$ such that
the vectors $J_i(a)$ are a basis of $T_{\gamma(a)}P$ and such that 
$J'_i(a)=-{\mathcal S}^P_{\dot\gamma(a)}(J_i(a))$. 
For $i=k+1,\ldots,m-1$, choose Jacobi fields $J_i$ such
that $J_i(a)=0$ and the vectors $J'_i(a)$ form a basis of $T_{\gamma(a)}P^\perp\cap\dot\gamma(a)^\perp$.
If $\gamma$ is lightlike choose
$J'_{m-1}(a)=\dot\gamma(a)$. Then, the $J_i$'s  form a basis of the space of $P$-Jacobi
fields orthogonal to $\gamma$. Now, we can write $V=\sum_{i=1}^{m-1}f_iJ_i$,
for piecewise smooth functions $f_i$. 

For, define $\bar J_i=J_i$ for
$i=1,\ldots,k$ and $\bar J_i(t)=J_i(t)/(t-a)$, $\bar J_i(a)=J'_i(a)$, for $i=k+1,\ldots,m-1$. The
absence of $P$-focal points along $\gamma$ and the fact that, under the hypothesis
that $\dot\gamma(a)\in T_{\gamma(a)}P^\perp\setminus T_{\gamma(a)}P$,
$T_{\gamma(a)}\M=T_{\gamma(a)}P\oplus T_{\gamma(a)}P^\perp$, imply that the vectors
$\bar J_i(t)$ are a basis for
$\dot\gamma(t)^\perp$ for $t\in[a,b]$.

Now, we have $J=\sum_{i=1}^{m-1}c_iJ_i$, where $c_i=f_i(b)$. The desired 
inequality follows 
directly from the Lemma~\ref{thm:lema} (equality implies that all $f_i$ are constant, 
except for $f_{m-1}$, in the lightlike
case).
\end{proof}
We give the
following definition:
\begin{defpartnorm}\label{thm:defpartnorm}
A partition $a=t_0<t_1<\ldots<t_N=b$ of $[a,b]$ is said to be {\em normal\/}
if the following conditions are satisfied:
\begin{itemize}
%\item[(a)] $\gamma(t_i)$ is not a $P$-focal point along $\gamma$ for
%all $i\ge1$;
\item[(a)] for all $i\ge1$ and all $t\in\,]t_i,t_{i+1}]$, the point
$\gamma(t)$ is not conjugate to $\gamma(t_i)$ along~$\gamma$;
\item[(b)] for all $t\in\,]t_0,t_1]$, the point $\gamma(t)$ is not $P$-focal
along $\gamma$.
\end{itemize}
\end{defpartnorm}
If $\gamma$ is either Riemannian or causal Lorentzian and $\dot\gamma(a)\in T_{\gamma(a)}P^\perp
\setminus T_{\gamma(a)}P$, since
the set of
$P$-focal points along $\gamma$ is finite, it is easy to see that
there exists $\delta>0$ such that every partition $t_0,\ldots,t_N$ of $[a,b]$
with $t_{i+1}-t_i\le\delta$ for all $i$ is normal. Namely, 
in order to (b) be satisfied, one can take $\delta$ to be the Lebesgue number
of a covering of $\gamma$ by {\em totally normal neighborhoods\/} 
(see~Ref.~\cite{dC}).

Given a normal partition, we define the following two subspaces 
of ${\mathcal H}^P$:
\begin{equation}\label{eq:defHoHj}
\begin{split}
&{\mathcal H}^P_0=\Big\{V\in{\mathcal H}^P:V(t_i)=0,\ \forall\,i\ge1\Big\};\\
&{\mathcal H}^P_J\!=\!\Big\{V\in{\mathcal H}^P\!:V\big\vert_{[t_i,t_{i+1}]}\
\text{is Jacobi}\ 
\forall\,i\ge1,\ \text{and}\ V\big\vert_{[t_0,t_1]}\ \text{is $P$-Jacobi}\Big\}.
\end{split}
\end{equation}
Observe that there exists an isomorphism:
\begin{equation}\label{eq:iso}
\phi:{\mathcal H}^P_J\longmapsto\bigoplus_{i=1}^{N-1}\dot\gamma(t_i)^\perp
\end{equation}
given by setting $\phi(V)=\big(V(t_1),V(t_2),\ldots, V(t_{N-1})\big)$.
Namely, since $\gamma(t_i)$ and $\gamma(t_{i+1})$ are non conjugate for
$i\ge1$, then $V\big\vert_{[t_i,t_{i+1}]}$ is uniquely determined
by the boundary values $V(t_i)$ and $V(t_{i+1})$; moreover, since
$\gamma(t_1)$ is not $P$-focal, then $V\big\vert_{[t_0,t_1]}$ is uniquely
determined by the value $V(t_1)$.

This shows that ${\mathcal H}^P_0\cap{\mathcal H}^P_J=\{0\}$ and 
that ${\mathcal H}^P_0+{\mathcal H}^P_J={\mathcal H}^P$, hence we have:
\begin{equation}\label{eq:somma}
{\mathcal H}^P_0\oplus{\mathcal H}^P_J={\mathcal H}^P.
\end{equation}
We are ready to prove the Morse Index Theorem for Riemannian
or causal Lorentz\-ian geodesics
with variable initial point:
\begin{indexth}\label{thm:indexth}
Let $(\M,g)$ be either Riemannian or Lorentzian, $P$ a smooth
submanifold of $\M$ and $\gamma:[a,b]\longmapsto\M$ a geodesic (causal, if $(M,g)$
is Lorentzian)
with $\gamma(a)\in P$ and $\dot\gamma(a)\in T_{\gamma(a)}P^\perp\setminus
T_{\gamma(a)}P$. Then, 
$\ind(I^P)=\sum\limits_{t_0\in\,]a,b[}\mu^P(t_0)<+\infty$.
\end{indexth}
\begin{proof}
For $[\alpha,\beta]\subset\,[\,a,b]$,
let $I_{[\alpha,\beta]}$ be the bilinear form (\ref{eq:bilform}) 
for the restricted geodesic $\gamma\vert_{[\alpha,\beta]}$
(omitting the term involving ${\mathcal S}^P_{\dot\gamma(a)}$); 
if $\alpha=a$, then we set $I^P_{[\alpha,\beta]}$ to be just the bilinear
form (\ref{eq:bilform}) for the restricted geodesic $\gamma\vert_{[\alpha,\beta]}$.
For $t\in]a,b]$ let's write $i(t)=\ind(I^P_{[a,t]})$; observe that $i(b)=\ind(I^P)$.
The function $i:[a,b]\longmapsto\N$ is non decreasing
(if $t<s$ we can regard $I^P_{[a,t]}$ as a restriction of
$I^P_{[a,s]}$, by extending vector fields on $[a,t]$ to $[a,s]$ defining 
them to be zero on $[t,s]$).

We show that $i(t)$ is piecewise constant and left-continuous on $[a,b]$,
and that $i(t^+)-i(t^-)=\mu^P(t)$ for all $t\in]a,b[$.

Let $t\in\,]a,b]$ be fixed and choose a normal partition
$t_0,\ldots,t_N$ of $[a,b]$ such that $t\in\,]t_i,t_{i+1}[$ for
some $i\ge1$ (we allow $t=t_{i+1}$ if $t=b$ and we set $i=N-1$). Let's denote by
${\mathcal H}^P_J([a,t])$
and ${\mathcal H}^P_0([a,t])$ the spaces defined in (\ref{eq:defHoHj}), 
replacing the  interval $[a,b]$ by $[a,t]$ (and using the normal partition
$t_0,\ldots,t_i,t$ of $[a,t]$).

\noindent\  
We observe that the direct sum (\ref{eq:somma}) (for the interval $[a,t]$) is 
$I^P_{[a,t]}$-orthogonal,
i.e., $I^P_{[a,t]}(V_0,V_J)=0$ for all $V_0\in{\mathcal H}^P_0([a,t])$ and 
$V_J\in{\mathcal H}^P_J([a,t])$,
which follows directly from (\ref{eq:bilformint}).

Next, we claim that $I^P_{[a,t]}\big\vert_{{\mathcal H}^P_0([a,t])}\ge0$. 
To check this, just observe that for $V\in{\mathcal H}^P_0([a,t])$ we have:
\[I^P_{[a,t]}(V,V)=I^P_{[t_0,t_1]}(V,V)+\sum_{j=1}^{i-1}I_{[t_j,t_{j+1}]}(V,V)+
I_{[t_i,t]}(V,V).\]
The claim now follows from
Corollary~\ref{thm:corJ}, by taking the Jacobi field $J=0$.

It follows that $i(t)=\ind(I^P_{[a,t]})=\ind(I^P_{[a,t]},{\mathcal H}^P_J([a,t]))$;
Observe that as in
(\ref{eq:iso}) the space ${\mathcal H}^P_J([a,t])$ is isomorphic to 
the space ${\mathcal H}_*$ defined by:
\[{\mathcal H}_*=\bigoplus_{j=1}^{i}\dot\gamma(t_j)^\perp,\]
and we'll call this isomorphism $\phi_t:{\mathcal
H}^P_J([a,t])\longmapsto{\mathcal H}_*$.

If $s\in]a,b]$ is sufficiently close to $t$ or, more precisely, if
$s\in]t_i,t_{i+1}]$, the arguments above can be repeated by replacing
$t$ with $s$ (observe, in particular, that the space ${\mathcal H}_*$ obtained will be
the same). We can use the isomorphism $\phi_s$ between ${\mathcal H}^P_J([a,s])$ and
${\mathcal H}_*$ to define a symmetric bilinear form $I_s$ on ${\mathcal H}_*$
corresponding to $I^P_{[a,s]}$. Clearly $i(s)=\ind(I_s)$.

We have now a one parameter family $I_s$ of symmetric bilinear forms on the (fixed) 
finite dimensional space ${\mathcal H}_*$ and it's not difficult to 
see that $I_s$ depends continuously
(actually, smoothly) on $s$.\footnote{To prove this fact, one uses equation (\ref{eq:bilformint}) 
to write a  expression for $I^P$ on piecewise Jacobi fields and  observes that the 
integral vanishes. Thus, formula~(\ref{eq:bilformint}) reduces to a  finite sum, and 
the conclusion follows from the theorem on smooth dependence on the initial
data for the solutions of the Jacobi differential equation.}

Let's consider the decomposition ${\mathcal H}_*={\mathcal H}_*^+
\oplus{\mathcal H}_*^-\oplus{\mathcal H}_*^0$, where $I_t$ is positive 
(respectively, negative) definite on
${\mathcal H}_*^+$ (respectively, ${\mathcal H}_*^-$) and ${\mathcal H}_*^0$ 
is the kernel of $I_t$. We can
also assume that this decomposition is $I_t$-orthogonal 
(this is just the Sylvester inertia Theorem). The
dimension of ${\mathcal H}_*^-$ is $i(t)$.

Since the decomposition ${\mathcal H}^P_0([a,t])\oplus{\mathcal H}^P_J([a,t])$
is orthogonal with respect to $I^P_{[a,t]}$, we know that the kernel of the
restriction of $I^P_{[a,t]}$ to ${\mathcal H}^P_J([a,t])$ (which corresponds
to ${\mathcal H}_*^0$ by the isomorphism $\phi_t$) is just the intersection of
${\mathcal H}^P_J([a,t])$ and the kernel of $I^P_{[a,t]}$, the last one being
given by Lemma \ref{thm:kernlema}. Observe that ${\mathcal
J}_0\subset{\mathcal H}^P_J([a,t])$ and denote by ${\mathcal J}_*$ the
subspace of ${\mathcal H}_*$ which corresponds to ${\mathcal J}_0$, i.e.,
${\mathcal J}_*=\phi_t({\mathcal J}_0)$. In the lightlike Lorentzian case,
write also ${\mathcal N}_*=\phi_t({\mathcal N}\cap{\mathcal H}^P_J([a,t]))$.

Observe that ${\mathcal N}_*$ is just the set of $i$-tuples of vectors which
are parallel to $\dot\gamma$, so that ${\mathcal N}_*$ doesn't change if we
replace $t$ by $s$ in its definition, and
therefore ${\mathcal N}_*$ is also contained in the kernel of $I_s$. 

We see now that ${\mathcal
H}_*^0={\mathcal J}_*$, except for the lightlike Lorentzian case where ${\mathcal H}_*^0={\mathcal J}_*\oplus{\mathcal N}_*$. The dimension of ${\mathcal J}_*$ is just the multiplicity $\mu^P(t)$ of $\gamma(t)$ as a $P$-focal point.

By the continuous dependence of $I_s$ on $s$ we see that for $\epsilon>0$
sufficiently small and $s\in]t-\epsilon,t+\epsilon[$, $I_s$ is negative
definite on ${\mathcal H}_*^-$ so that $i(s)\geq i(t)$. For
$s\in]t-\epsilon,t]$ we have also $i(s)\leq i(t)$ so that $i(s)=i(t)$, i.e.,
$i$ is constant on $]t-\epsilon,t]$. This finishes the proof that $i$ is left
continuous. From now on we suppose $t<b$.

The same continuity argument show that for some $\epsilon>0$, we have that
$I_s$ is positive definite on ${\mathcal H}_*^+$ for $s\in[t,t+\epsilon[$ (and
positive semi-definite on ${\mathcal H}_*\oplus{\mathcal N}_*$ for $\gamma$ lightlike), so that $i(s)$ is bounded above by the codimension of
${\mathcal H}_*^+$ (or ${\mathcal H}_*^+\oplus{\mathcal N}_*$, respectively). If $\gamma(t)$ is not a $P$-focal point this codimension equals $i(t)$ so that $i(s)=i(t)$ for $s\in]t-\epsilon,t+\epsilon[$.

Finally, if $\gamma(t)$ is a $P$-focal point, by the above argument 
we only obtain the inequality $i(s)\leq
i(t)+\mu^P_\gamma(t)$. We'll show below that for $s\in]t,t_{i+1}]$ and for
$V=(v_1,\ldots,v_i)\in{\mathcal H}_*$ we have $I_s(V,V)\leq I_t(V,V)$, the
inequality being strict if $v_i\not=0$ (or if $v_i$ is not parallel to $\dot\gamma$,
in case $\gamma$ is lightlike). But this hypothesis on $v_i$ holds
if $V\in{\mathcal J}_*$ and $V\not=0$, observing that the corresponding vector
field $\phi_t^{-1}(V)$ on ${\mathcal H}^P_J([a,t])$ is an unbroken Jacobi field. 
We conclude then that $I_s(V,V)<0$ for nonzero $V\in{\mathcal J}_*$ and hence 
for all nonzero $V\in{\mathcal H}_*^-\oplus{\mathcal J}_*$, which implies 
that $I_s$ is negative definite on this space and $i(s)\geq i(t)+\mu^P(t)$.

We are now left with the proof of the inequality $I_s(V,V)\leq I_t(V,V)$. Towards this
goal, let 
$V_1\in{\mathcal H}^P_J([a,t])$ and $V_2\in{\mathcal H}^P_J([a,s])$ be the
vector fields corresponding to $V\in{\mathcal H}_*$, i.e.,
$V_1=\phi_t^{-1}(V)$ and $V_2=\phi_s^{-1}(V)$. Extend $V_1$ to  zero on $[t,s]$. 
Then, $I_t(V,V)=I^P_{[a,s]}(V_1,V_1)$ and $I_s(V,V)=I^P_{[a,s]}(V_2,V_2)$. The
vector fields $V_1$ and $V_2$ differ at the most in the interval $[t_i,s]$. 
The restriction of $V_1$ to $[t_i,t]$ is the only Jacobi field such 
that $V_1(t_i)=v_i$ and $V_1(t)=0$,
while the restriction of $V_2$ to $[t_i,s]$ is the only Jacobi field 
such that $V_2(t_i)=v_i$ and
$V_2(s)=0$. We have:
\[I_t(V,V)-I_s(V,V)=I_{[t_i,s]}(V_1,V_1)-I_{[t_i,s]}(V_2,V_2).\]
We now apply Corollary \ref{thm:corJ} to the geodesic
$\gamma\vert_{[t_i,s]}$ (with starting and ending points interchanged), for
the Jacobi field $V_2$, vector field $V_1$ and submanifold equal to the point
$\{\gamma(s)\}$. For the strict inequality we need the hypothesis that
$v_i\not=0$ (respectively, $v_i$ not parallel to $\dot\gamma$, in the
lightlike Lorentzian case), since this implies that $V_1$ is not Jacobi in
$[t_i,s]$ (respectively, does not differ from a Jacobi field by a multiple of
$\dot\gamma$, in the lightlike Lorentzian case). This concludes the proof.
\end{proof}
\begin{light}\label{thm:light}
If $(M,g)$ is Lorentzian, then the case that $\dot\gamma(a)\in T_{\gamma(a)}P\cap
T_{\gamma(a)}P^\perp$ may happen only when $\gamma$ is lightlike and
$P$ is a degenerate submanifold at $\gamma(a)$, i.e., the restriction of 
$g$ to $T_{\gamma(a)}P$ is degenerate. Observe that in this case the thesis of
Theorem~\ref{thm:indexth} is clearly false. For instance, if $\M=\R^2$ and $g$
is the flat Minkowski metric ${\rm d}x^2-{\rm d}t^2$, $P$ is the diagonal
$x=t$ and $\gamma$ is any segment contained in $P$, then every point of
$\gamma$ is $P$-focal.
\end{light}

We now want to extend the Morse Index Theorem to the case of two
variable endpoints. To this end, we now assume that $P$ and $Q$
are smooth submanifolds of $\M$, and that $\gamma:[a,b]\longmapsto\M$
is a geodesic with $\gamma(a)\in P$, $\dot\gamma(a)\in T_{\gamma(a)}P^\perp$,
$\gamma(b)\in Q$ and $\dot\gamma(b)\in T_{\gamma(b)}Q^\perp$.

We denote by ${\mathcal H}^{(P,Q)}$ the vector space of all piecewise smooth
vector fields $V$ along $\gamma$, with $g(V,\dot\gamma)\equiv0$, $V(a)\in T_{\gamma(a)}P$
and $V(b)\in T_{\gamma(b)}Q$. Moreover, we will consider the symmetric bilinear
form $I^{(P,Q)}$ on ${\mathcal H}^{(P,Q)}$, given by:
\begin{equation}\label{eq:defIPQ}I^{(P,Q)}(V,W)=I^P(V,W)+{\mathcal
S}^Q_{\dot\gamma(b)}(V(b),W(b)).
\end{equation} 
%A $Q$-Jacobi field along $\gamma$ is defined to be a Jacobi
%field $J$ along $\gamma$, with $J(b)\in T_{\gamma(b)}Q$ and such that (\ref{eq:forgot})
%is satisfied when $P$ is replaced by $Q$ and $a$ is replaced by $b$. 
%
%Observe that a vector field $V\in {\mathcal H}^{(P,Q)}$ is in the kernel
%of $I^{(P,Q)}$ if and only if $V$ is at the same time a $P$-Jacobi field
%and a $Q$-Jacobi field. 

Let ${\mathcal J}^Q$ denote the subspace of ${\mathcal H}^{(P,Q)}$ consisting
of $P$-Jacobi fields, and let ${\mathcal A}$ be the symmetric bilinear
form on $\mathcal J^Q$ obtained by the restriction of $I^{(P,Q)}$.
Then, it is easily computed from \eqref{eq:bilform} using integration by parts:
\[{\mathcal A}(J_1,J_2)={\mathcal S}^Q_{\dot\gamma(b)}(J_1(b),J_2(b))+
g(J_1'(b),J_2(b)),\quad J_1,J_2\in{\mathcal J}^Q.\]
Moreover, for $t\in[a,b]$, we introduce the space ${\mathcal J}[t]$:
\[{\mathcal J}[t]=\Big\{J(t):J\ \text{is $P$-Jacobi}\ \Big\}\subset T_{\gamma(t)}\M;\]
observe that, for $t\in\,]a,b]$, $\gamma(t)$ is not $P$-focal if and only if  ${\mathcal
J}[t]=T_{\gamma(t)}\M$.

We can now state and prove the following extension
of the Morse Index Theorem for geodesics between
submanifolds:
\begin{twosubmanifolds}\label{thm:twosubmanifolds} 
Let $(\M,g)$ be a semi-Riemannian manifold, $P,Q$ submanifolds of $\M$ and
$\gamma:[a,b]\longmapsto\M$ be a geodesic such that $\gamma(a)\in P$,
$\dot\gamma(a)\in T_{\gamma(a)}P^\perp$, $\gamma(b)\in Q$ and
$\dot\gamma(b)\in T_{\gamma(b)}Q^\perp$. 
Assume that ${\mathcal J}[b]\supset T_{\gamma(b)}Q$.
%$\gamma(b)$ is not a $P$-focal point along $\gamma$.
Let ${\mathcal V}$ be a subspace
of ${\mathcal H}^{(P,Q)}$ that contains the space ${\mathcal J}^Q$ of
$P$-Jacobi fields along $\gamma$ in ${\mathcal H}^{(P,Q)}$. Then,
$\ind(I^{(P,Q)},{\mathcal V})=\ind(I^P,{\mathcal H}^P\cap{\mathcal V})+\ind({\mathcal A},
{\mathcal J})$.
\end{twosubmanifolds}
\begin{proof}
The space ${\mathcal H}^P$ is given by the subspace of ${\mathcal H}^{(P,Q)}$
consisting of those vector fields $V$ such that $V(b)=0$; moreover, the
restriction of $I^{(P,Q)}$ to ${\mathcal H}^P$ is precisely $I^P$. 
Defining $\mathcal J_0$ as in formula \eqref{eq:defN}, 
let ${\mathcal J}_1$ be any subspace of
${\mathcal J}^Q$ such that ${\mathcal J}^Q={\mathcal J}_0\oplus{\mathcal J}_1$.
Clearly, ${\mathcal H}^{(P,Q)}={\mathcal H}^P\oplus{\mathcal J}_1$, because
${\mathcal J}[b]\supset T_{\gamma(b)}Q$;
moreover, from (\ref{eq:defIPQ}) it follows immediately that
this decomposition is $I^{(P,Q)}$-orthogonal, i.e., $I^{(P,Q)}(V,J)=0$
for all $V\in {\mathcal H}^P$ and all $J\in{\mathcal J}_1$.
Since ${\mathcal V}$ contains ${\mathcal J}_1$, then ${\mathcal V}=({\mathcal V}\cap
{\mathcal H}^P)\oplus{\mathcal J}_1$.
Hence, $\ind(I^{(P,Q)},{\mathcal V})=\ind(I^P,{\mathcal H}^P\cap{\mathcal V})
+\ind({\mathcal A}, {\mathcal J}_1)$. To conclude the proof, we simply observe that
$\ind({\mathcal A}, {\mathcal J}_1)=\ind({\mathcal A}, {\mathcal J})$,
because ${\mathcal J}_0\subset{\rm Ker}({\mathcal A})$.
\end{proof}
\begin{h1}\label{thm:h1}
One can consider suitable Hilbert space completions $\bar{\mathcal H}^P$
and $\bar{\mathcal H}^{(P,Q)}$ of the spaces
${\mathcal H}^P$ and ${\mathcal H}^{(P,Q)}$ with respect to an $H^1$-Sobolev
norm. Then, the bilinear forms $I^P$ and $I^{(P,Q)}$ extend uniquely
to bounded symmetric bilinear forms on these Hilbert spaces.
Observe that a bounded symmetric bilinear form 
on a Hilbert space and its restriction to any dense subspace
have the same index. Using a Hilbert space approach, Theorem~\ref{thm:indexth}
can be proven alternatively by means of the spectral theory for compact
self-adjoint operators (see \cite[Theorem~5.9.3]{Ma} for an idea of such a proof).
\end{h1}
\begin{oldresults}\label{thm:oldresults}
If $(\M,g)$ is Riemannian and ${\mathcal V}={\mathcal H}^{(P,Q)}$,
then Theorems~\ref{thm:indexth} and ~\ref{thm:twosubmanifolds} give as a particular case the
Index Theorem of~\cite[p.\ 342]{Kal} and the older versions of the Morse Index
Theorem presented in \cite{A,Bol}. 
In~\cite{EK} it was briefly mentioned the fact that 
results analogous to the Riemannian case could apply to
the Lorentzian timelike case.
As to the lightlike case, in References~\cite{BE, BEE, EK} the authors 
consider the index of $I^{(P,Q)}$ in the quotient space ${\mathcal H}^{(P,Q)}/{\mathcal N}$
(recall formula~(\ref{eq:defN})); in this situation, ${\mathcal N}$ is contained
in the kernel of $I^{(P,Q)}$. By a simple linear algebra argument
one proves that the index of a bilinear form in a quotient space by a subspace
of its kernel is the same as the index of the form in the original space.
Hence, Theorems~\ref{thm:indexth} and \ref{thm:twosubmanifolds} generalize
the results of \cite{BE, BEE, EK}.
\end{oldresults}
\begin{stationary}\label{thm:stationary}
The result of Theorem~\ref{thm:twosubmanifolds} becomes significant
when the subspace ${\mathcal V}$ of ${\mathcal H}^{(P,Q)}$
is chosen in such a way that $\ind(I^P,{\mathcal
H}^P\cap{\mathcal V})$ is finite; observe that $\ind({\mathcal A},{\mathcal J})$ 
is always finite. If one considers geodesics in semi-Riemannian manifolds with metric
of index greater or equal to $2$, or spacelike geodesics in Lorentzian manifolds,
then $\ind({I^P,\mathcal H}^P)$ is in general infinite (see Refs.~\cite{BGM, Hel} for 
further results in this direction). Nevertheless, the restriction to suitable subspaces
may yield the finiteness of the index, and, possibly, weaker versions of the Morse
Index Theorem may apply. 
For instance (see Ref.~\cite{GP}), let's consider the
case  of a stationary Lorentzian manifold $(\M,g)$, i.e., a Lorentzian manifold endowed with a
timelike Killing vector field $Y$. Let $\gamma$ be a spacelike geodesic; we consider for
simplicity the case that the initial manifold $P$ reduces to a point.
The Killing vector field $Y$ induces the conservation law $g(\dot\gamma,Y)\equiv C_\gamma$
for all geodesic $\gamma$; then, one can consider only variations $\gamma_s$
of $\gamma$ such that $g(\dot\gamma_s,Y)\equiv C_s$, and the corresponding
variational field $V=\ddso\gamma_s$ belongs to the space:
\[{\mathcal V}=\Big\{V:\exists\,C_V\in\R\ \text{such that}\ g(V',Y)-g(V,Y')\equiv C_V\Big\}.\]
It is a simple observation that $\mathcal V$ contains all the Jacobi
fields along $\gamma$; moreover, using the Sobolev Embedding Theorem
one proves that the bilinear form $I^P$ is given by a self-adjoint
operator $T$ on the closure of $\mathcal V\cap {\mathcal H}^P$ in a 
suitable Sobolev space completion  of ${\mathcal H}^P$,
where $T$ is a {\em compact perturbation\/} of the identity (see Ref.~\cite{Ma}). 
Thus,  $\mathcal V\cap {\mathcal H}^{(P,Q)}$ satisfies the hypothesis of 
Theorem~\ref{thm:twosubmanifolds}
and it is such that $\ind(I^P,{\mathcal V}\cap{\mathcal H}^P)$ is finite. The question of whether
such index equals the geometrical index of $\gamma$ remains still unanswered.\end{stationary}
\end{section}
%%%%%%%%%%%%%%%%%%%%%%%%%%%%%%%%%%
%%%%%%%%%%%%%%%%%%%%%%%%%%%%%%%%%%

%%%%%%%%%%%%%%%%%%%%%%%%%%%%%%%%%%%%%%%%%%%%%%%%%%%%%%%%%%%%%%%%%%%%%%%%%%
%%%%%%%%%%%%%%%%%%%%%%%%%%%%%%%%%%%%%%%%%%%%%%%%%%%%%%%%%%%%%%%%%%%%%%%%%%

\end{document}